\documentclass{amsart}
\usepackage{amssymb,amsmath,theorem,euscript}

\input{epsf}

\newcounter{punct} \def\punct{\addtocounter{punct}{1}{ \arabic{punct}}}

\def\sm{\smallskip}

\begin{document}

\def\wh{\widehat}

 \def\ov{\overline}
\def\wt{\widetilde}
 \newcommand{\rk}{\mathop {\mathrm {rk}}\nolimits}
\newcommand{\Aut}{\mathop {\mathrm {Aut}}\nolimits}
\newcommand{\Out}{\mathop {\mathrm {Out}}\nolimits}
 \newcommand{\tr}{\mathop {\mathrm {tr}}\nolimits}
  \newcommand{\diag}{\mathop {\mathrm {diag}}\nolimits}
  \newcommand{\supp}{\mathop {\mathrm {supp}}\nolimits}
  \newcommand{\indef}{\mathop {\mathrm {indef}}\nolimits}
  \newcommand{\dom}{\mathop {\mathrm {dom}}\nolimits}
  \newcommand{\im}{\mathop {\mathrm {im}}\nolimits}
 
\renewcommand{\Re}{\mathop {\mathrm {Re}}\nolimits}

\def\Br{\mathrm {Br}}

\def\SL{\mathrm {SL}}
\def\Diag{\mathrm {Diag}}
\def\SU{\mathrm {SU}}
\def\GL{\mathrm {GL}}
\def\U{\mathrm U}
\def\OO{\mathrm O}
 \def\Sp{\mathrm {Sp}}
 \def\SO{\mathrm {SO}}
\def\SOS{\mathrm {SO}^*}
 \def\Diff{\mathrm{Diff}}
 \def\Vect{\mathfrak{Vect}}
\def\PGL{\mathrm {PGL}}
\def\PU{\mathrm {PU}}
\def\PSL{\mathrm {PSL}}
\def\Symp{\mathrm{Symp}}
\def\End{\mathrm{End}}
\def\Mor{\mathrm{Mor}}
\def\Aut{\mathrm{Aut}}
 \def\PB{\mathrm{PB}}
 \def\cA{\mathcal A}
\def\cB{\mathcal B}
\def\cC{\mathcal C}
\def\cD{\mathcal D}
\def\cE{\mathcal E}
\def\cF{\mathcal F}
\def\cG{\mathcal G}
\def\cH{\mathcal H}
\def\cJ{\mathcal J}
\def\cI{\mathcal I}
\def\cK{\mathcal K}
 \def\cL{\mathcal L}
\def\cM{\mathcal M}
\def\cN{\mathcal N}
 \def\cO{\mathcal O}
\def\cP{\mathcal P}
\def\cQ{\mathcal Q}
\def\cR{\mathcal R}
\def\cS{\mathcal S}
\def\cT{\mathcal T}
\def\cU{\mathcal U}
\def\cV{\mathcal V}
 \def\cW{\mathcal W}
\def\cX{\mathcal X}
 \def\cY{\mathcal Y}
 \def\cZ{\mathcal Z}
\def\0{{\ov 0}}
 \def\1{{\ov 1}}
 \def\frA{\mathfrak A}
 \def\frB{\mathfrak B}
\def\frC{\mathfrak C}
\def\frD{\mathfrak D}
\def\frE{\mathfrak E}
\def\frF{\mathfrak F}
\def\frG{\mathfrak G}
\def\frH{\mathfrak H}
\def\frI{\mathfrak I}
 \def\frJ{\mathfrak J}
 \def\frK{\mathfrak K}
 \def\frL{\mathfrak L}
\def\frM{\mathfrak M}
 \def\frN{\mathfrak N} \def\frO{\mathfrak O} \def\frP{\mathfrak P} \def\frQ{\mathfrak Q} \def\frR{\mathfrak R}
 \def\frS{\mathfrak S} \def\frT{\mathfrak T} \def\frU{\mathfrak U} \def\frV{\mathfrak V} \def\frW{\mathfrak W}
 \def\frX{\mathfrak X} \def\frY{\mathfrak Y} \def\frZ{\mathfrak Z} \def\fra{\mathfrak a} \def\frb{\mathfrak b}
 \def\frc{\mathfrak c} \def\frd{\mathfrak d} \def\fre{\mathfrak e} \def\frf{\mathfrak f} \def\frg{\mathfrak g}
 \def\frh{\mathfrak h} \def\fri{\mathfrak i} \def\frj{\mathfrak j} \def\frk{\mathfrak k} \def\frl{\mathfrak l}
 \def\frm{\mathfrak m} \def\frn{\mathfrak n} \def\fro{\mathfrak o} \def\frp{\mathfrak p} \def\frq{\mathfrak q}
 \def\frr{\mathfrak r} \def\frs{\mathfrak s} \def\frt{\mathfrak t} \def\fru{\mathfrak u} \def\frv{\mathfrak v}
 \def\frw{\mathfrak w} \def\frx{\mathfrak x} \def\fry{\mathfrak y} \def\frz{\mathfrak z} \def\frsp{\mathfrak{sp}}
 \def\bfa{\mathbf a} \def\bfb{\mathbf b} \def\bfc{\mathbf c} \def\bfd{\mathbf d} \def\bfe{\mathbf e} \def\bff{\mathbf f}
 \def\bfg{\mathbf g} \def\bfh{\mathbf h} \def\bfi{\mathbf i} \def\bfj{\mathbf j} \def\bfk{\mathbf k} \def\bfl{\mathbf l}
 \def\bfm{\mathbf m} \def\bfn{\mathbf n} \def\bfo{\mathbf o} \def\bfp{\mathbf p} \def\bfq{\mathbf q} \def\bfr{\mathbf r}
 \def\bfs{\mathbf s} \def\bft{\mathbf t} \def\bfu{\mathbf u} \def\bfv{\mathbf v} \def\bfw{\mathbf w} \def\bfx{\mathbf x}
 \def\bfy{\mathbf y} \def\bfz{\mathbf z} \def\bfA{\mathbf A} \def\bfB{\mathbf B} \def\bfC{\mathbf C} \def\bfD{\mathbf D}
 \def\bfE{\mathbf E} \def\bfF{\mathbf F} \def\bfG{\mathbf G} \def\bfH{\mathbf H} \def\bfI{\mathbf I} \def\bfJ{\mathbf J}
 \def\bfK{\mathbf K} \def\bfL{\mathbf L} \def\bfM{\mathbf M} \def\bfN{\mathbf N} \def\bfO{\mathbf O} \def\bfP{\mathbf P}
 \def\bfQ{\mathbf Q} \def\bfR{\mathbf R} \def\bfS{\mathbf S} \def\bfT{\mathbf T} \def\bfU{\mathbf U} \def\bfV{\mathbf V}
 \def\bfW{\mathbf W} \def\bfX{\mathbf X} \def\bfY{\mathbf Y} \def\bfZ{\mathbf Z} \def\bfw{\mathbf w}
 \def\R {{\mathbb R }} \def\C {{\mathbb C }} \def\Z{{\mathbb Z}} \def\H{{\mathbb H}} \def\K{{\mathbb K}}
 \def\N{{\mathbb N}} \def\Q{{\mathbb Q}} \def\A{{\mathbb A}} \def\T{\mathbb T} \def\P{\mathbb P} \def\G{\mathbb G}
 \def\bbA{\mathbb A} \def\bbB{\mathbb B} \def\bbD{\mathbb D} \def\bbE{\mathbb E} \def\bbF{\mathbb F} \def\bbG{\mathbb G}
 \def\bbI{\mathbb I} \def\bbJ{\mathbb J} \def\bbK{\mathbb K} \def\bbL{\mathbb L} \def\bbM{\mathbb M} \def\bbN{\mathbb N} \def\bbO{\mathbb O}
 \def\bbP{\mathbb P} \def\bbQ{\mathbb Q} \def\bbS{\mathbb S} \def\bbT{\mathbb T} \def\bbU{\mathbb U} \def\bbV{\mathbb V}
 \def\bbW{\mathbb W} \def\bbX{\mathbb X} \def\bbY{\mathbb Y} \def\kappa{\varkappa} \def\epsilon{\varepsilon}
 \def\phi{\varphi} \def\le{\leqslant} \def\ge{\geqslant}

\def\UU{\bbU}
\def\Mat{\mathrm{Mat}}
\def\tto{\rightrightarrows}

\def\Gr{\mathrm{Gr}}

\def\graph{\mathrm{graph}}

\def\O{\mathrm{O}}

\def\la{\langle}
\def\ra{\rangle}

\def\B{\mathrm B}
\def\Int{\mathrm{Int}}
\def\LGr{\mathrm{LGr}}


\def\I{\mathbb I}
\def\M{\mathbb M}
\def\T{\mathbb T}

\def\Lat{\mathrm{Lat}}
\def\LLat{\mathrm{LLat}} 
\def\Mod{\mathrm{Mod}}
\def\LMod{\mathrm{LMod}}
\def\Naz{\mathrm{Naz}}
\def\naz{\mathrm{naz}}
\def\bNaz{\mathbf{Naz}}
\def\AMod{\mathrm{AMod}}
\def\ALat{\mathrm{ALat}}
\def\MAT{\mathrm{MAT}}
\def\Mar{\mathrm{Mar}}

\def\Ver{\mathrm{Vert}}
\def\Bd{\mathrm{Bd}}
\def\We{\mathrm{We}}
\def\Heis{\mathrm{Heis}}
\def\Pol{\mathrm{Pol}}
\def\Ams{\mathrm{Ams}}
\def\Herm{\mathrm{Herm}}

 \def\kos{/\!\!/}
  \newcommand{\sgn}{\mathop {\mathrm {sgn}}\nolimits}
  
  \renewcommand{\Re}{\mathop {\mathrm {Re}}\nolimits}
\renewcommand{\Im}{\mathop {\mathrm {Im}}\nolimits}

\def\Gr{\mathrm{Gr}}

\begin{center}
 \Large\bf

Operational calculus for Fourier transform on the group $\GL(2,\R)$

\bigskip

\large \sc

\sm

Yury A. Neretin%
 \footnote{Supported by the grant FWF, P28421.}

\end{center}

\vspace{22pt}

{\small 
Consider the  Fourier transform on the group $\GL(2,\R)$ of real $2\times 2$-matrices.
We show that Fourier-images of polynomial differential
operators on $\GL(2,\R)$ are  differential-difference operators
with coefficients meromorphic in parameters of representations.
Expressions for operators contain shifts in imaginary direction with respect to the
integration contour in the Plancherel formula. We present explicit formulas
for images of partial derivations and multiplications by coordinates.}

\medskip

The  Plancherel formula for $\SL(2,\C)$ for obtained in by I.~M.~Gelfand and M.~I.~Nai\-mark
in 1947, for $\SL(2,\R)$ by Harich-Chandra in 1951, see \cite{GN}, \cite{Harish}. 
Later there were published many Plancherel formulas 
for Fourier transforms on non-commutative locally compact groups and homogeneous spaces 
(see a collection of references in \cite{Ner-after}). 
The spherical transform for Riemannian symmetric spaces and for Bruhat--Tits buildings
had numerous continuations in mathematics,
however a destiny of the majority of Plancherel formulas seems strange, they are heavy impressive results, 
which are difficult 
for usage. The corresponding integral transforms  have properties that immediately
follow from their definition.
Efforts of a further expansion
are faced with difficulties,
also there is a few number of explicit calculations of Fourier-images of functions.

An example of a nontrivial transformation of  differential operators 
for an $\SL(2,\R)$-related Fourier transform was done in \cite{Ner-imaginary}, 
images of certain first-order differential operators  are differential-difference operators
with  differentiations of the  order two in a space variable and shifts in imaginary direction in a parameter variable.
Several correspondences of the same type were obtained bu V.~F.~Molchanov and the author in
\cite{Mol1}--\cite{Mol3}, \cite{Ner-Gl}--\cite{Ner-jfaa}.
The present  paper contains a simple corollary of calculations of \cite{Ner-jfaa}, however 
our  statement (Theorems 1-2) is essentially stronger and easier for a further usage than the statement
of \cite{Ner-jfaa}.

\sm

{\bf \punct. The principal series of representations of the group $\GL(2,\R)$.}
Consider    the group $\GL(2,\R)$ of  real invertible matrices $\begin{pmatrix}
                                                                  a&b\\c&d
                                                                 \end{pmatrix}$
of order $2$.
For $\mu\in\C$ and $\epsilon\in\Z_2$ we define the function
$x^{\mu\kos  \epsilon}$ on $\R\setminus 0$ by
$$
x^{\mu\kos  \epsilon}:=|x|^\mu \sgn(x)^\epsilon.
$$

Denote 
$\Lambda:=\C\times \Z_2\times \C\times \Z_2$.
For each element $(\mu_1,\epsilon_1;\mu_2,\epsilon_2)$ of $\Lambda$ we define a representation  
$T_{\mu,\epsilon}$  of $\GL_2(\R)$ in the space of functions 
on $\R$ by 
\begin{multline*}
T_{\mu_1,\epsilon_1;\mu_2,\epsilon_2}\begin{pmatrix}
a&b\\
c&d
\end{pmatrix} \phi(t)=\\=
\phi\Bigl(\frac{b+t d}{a+t c}\Bigr)
\cdot(a+t c)^{-1+\mu_1-\mu_2\kos  \epsilon_1-\epsilon_2} \det\begin{pmatrix}
a&b\\
c&d
\end{pmatrix}^{1/2+\mu_2\kos  \epsilon_2}
.
\end{multline*}

This formula determines the (non-unitary) {\it principal series of representations of $\GL(2,\R)$}.
To be definite, we define a  representation $T_{\mu_1,\epsilon_1;\mu_2,\epsilon_2}$
in the space $C^{\infty}_{\mu_1-\mu_2,\epsilon_1-\epsilon_2}(\R)$ of smooth functions $f$ on $\R$
satisfying the additional condition on asymptotics at infinity:
the function 
$$
T_{\mu_1,\epsilon_1;\mu_2,\epsilon_2}\begin{pmatrix}
0&1\\
-1&0
\end{pmatrix} f(t)=f(-t^{-1}) (-t)^{-1+\mu_1-\mu_2\kos  \epsilon_1-\epsilon_2}
$$
must be smooth at 0. This space is invariant with respect to the whole group
$\GL(2,\R)$, see \cite{GGV}, Subsect. VII.1.2.

If $\mu_1-\mu_2-1\notin 2\Z+(\epsilon_1-\epsilon_2)$, then representations  $T_{\mu_1,\epsilon_1;\mu_2,\epsilon_2}$ and
$T_{\mu_2,\epsilon_2;\mu_1,\epsilon_1}$ are irreducible and equivalent 
(see, e.g., \cite{GGV}, Subsect. VII.2.2). The intertwining operator is given by 
$$
A_{\mu_1,\epsilon_1;\mu_2,\epsilon_2} f(t):=
\int_{\R} (t-s)^{-1-\mu_1+\mu_2\kos  \epsilon_1-\epsilon_2}
f(s)\,ds.
$$
The integral converges if $\Re(\mu_2-\mu_1)>0$ and is holomorphic in the parameters $\mu_1$, $\mu_2$.
It has a meromorphic continuation to the whole domain $\C\times \C$, for $\epsilon_1-\epsilon_2=0$
possible poles are on the complex $\mu_2-\mu_1=-1$, $-3$, \dots; for $\epsilon_1-\epsilon_2=1$
on $\mu_2-\mu_1=0$, 2, 4, \dots.

If $\mu_1=i \tau_1$, $\mu_2=i \tau_2\in i\R$, then a representation $T_{\mu_1,\epsilon_1;\mu_2,\epsilon_2}$ 
is unitary in $L^2(\R)$
(they are called  representations of the {\it unitary principal series}).

\sm

{\bf \punct. Fourier transform in the complex domain.} Consider the space
$C_0^\infty\bigl(\GL(2,\R)\bigr)$  of compactly supported smooth functions on $\GL(2,\R)$.
For $F\in C_0^\infty\bigl(\GL(2,\R)\bigr)$ we define its Fourier transform as 
a function sending $(\mu_1,\epsilon_1, \mu_2, \epsilon_2)$ to an operator
in $C^\infty_{\mu_1-\mu_2,\epsilon_1-\epsilon_2}(\R)$ given by
$$
T_{\mu_1,\epsilon_1;\mu_2,\epsilon_2}(F)  f(t)=
\int_{\GL_2(\R)} F_{\mu_1,\epsilon_1;\mu_2,\epsilon_2}(g)\, T(g)\, dg,
$$
where $dg$ denotes the Haar measure on $\GL_2(\R)$, its density with respect to
the standard Lebesgue  measure on the space of matrices is $(ad-bc)^{-2}$.

We denote the space of operator-valued functions on $\Lambda$
obtained in this way by $\cE$. A precise description of $\cE$  (a {\it Paley--Wiener theorem} for $\SL(2,\R)$)
was obtained 
in \cite{EM}, we only present some remarks.

Elements of $\cE$ satisfy a symmetry condition
$$
A_{\mu_1,\epsilon_1;\mu_2,\epsilon_2}\,T_{\mu_1,\epsilon_1;\mu_2,\epsilon_2}(F)
=
T_{\mu_2,\epsilon_2;\mu_1,\epsilon_1}(F)\,A_{\mu_1,\epsilon_1;\mu_2,\epsilon_2}.
$$

For any $F\in C_0^\infty\bigl(\GL(2,\R)\bigr)$ the operators $T_{\mu_1,\epsilon_1;\mu_2,\epsilon_2}(F)$
can be represented as integral operators
$$
T_{\mu_1,\epsilon_1;\mu_2,\epsilon_2}(F)f(t)=
\int_\R K_F(t,s| \mu_1,\epsilon_1;\mu_2,\epsilon_2)\, f(s)\, ds
,$$
where
 	\begin{multline*}
K_F(t,s|\mu_1,\epsilon_1;\mu_2,\epsilon_2)=\\=
\iiint\limits_{\R^3} F\bigl(
u - t v, s u - s t v - t w, v, s v + w \bigr)\,u^{-3/2+\mu_1\kos  \epsilon_1}
w^{-3/2+\mu_2\kos \epsilon_2}\,du\,dv\,dw.
\end{multline*}
It is easy to see  that  a support of the integrand is 
a bounded set having empty intersections with planes $u=0$, $w=0$.
Therefore 
the function $K_F$ is holomorphic in $\mu_1$, $\mu_2\in \C$ and smooth in
$s$, $t\in\R$. Also it satisfies some conditions on asymptotics at infinity 
(see a discussion of a similar case $\SL(2,\C)$ in \cite{GGV}, \S IV.1).

\sm

{\bf \punct. Immediate properties of the Fourier transform.}
Denote by $F_1*F_2$ the convolution on $\GL(2,\R)$. For any
$F_1$, $F_2\in C_0^\infty\bigl(\GL(2,\R)\bigr)$ we have
$$
T(F_1)\,T(F_2)=T(F_1*F_2).
$$
This is a general fact for locally compact 
groups with two-side invariant Haar measure, see e.g., \cite{Kir}, \S 10.

Next, the Lie algebra $\mathfrak{gl}(2)$
acts in $C^\infty\bigl(\GL(2,\R)\bigr)$ by right-invariant vector fields, we write formulas 
for the standard generators:
 \begin{align}
  \label{eq:er1}
e_{11}^r= - a \frac\partial{\partial a} - b \frac\partial{\partial b}
,\qquad
e_{12}^r= -c \frac\partial{\partial a} -  d \frac\partial{\partial b},\\
 \label{eq:er2}
e_{21}^r= -a \frac\partial{\partial c}  - b \frac\partial{\partial d},\qquad
e_{22}^r= -c \frac\partial{\partial c} - d \frac\partial{\partial d}.
\end{align}
Also, $\mathfrak{gl}(2)$
acts in $C^\infty\bigl(\GL(2,\R)\bigr)$ by left-invariant vector fields
\begin{align}
e_{11}^l= a \frac\partial{\partial a} + c \frac\partial{\partial c},\qquad
 \label{eq:el1}
e_{12}^l= a \frac\partial{\partial b} + c \frac\partial{\partial d},\\
 \label{eq:el2}
e_{21}^l= b \frac\partial{\partial a}  + d \frac\partial{\partial c} ,\qquad
e_{22}^l= b \frac\partial{\partial b}  + d \frac\partial{\partial d}.
\end{align}
Thus we get an action of the Lie algebra $\mathfrak{gl}(2)\oplus \mathfrak{gl}(2)$
in $C^\infty\bigl(\GL(2,\R)\bigr)$.
The corresponding operators in the Fourier-image are given by
 \begin{align}
E_{11}^r&= -t \frac{\partial}{\partial t}-(1/2-\mu_1),
\qquad &E_{12}^r&=\frac{\partial}{\partial t},
\label{eq:Er1}
\\
E_{21}^r&=-t^2 \frac{\partial}{\partial t}+ (-1+\mu_1-\mu_2)t,
\qquad
&E_{22}^r&= t \frac{\partial}{\partial t}+(1/2 + \mu_2),
\label{eq:Er2}
\end{align}
and
 \begin{align}
E_{11}^l&=-s\frac\partial{\partial s}- (1/2+\mu_1),
\qquad
&E_{12}^l&=\frac\partial{\partial s},
\label{eq:El1}
\\
E_{21}^l&=-s^2\frac\partial{\partial s}+(-1-\mu_1+\mu_2)s,
\qquad
& E_{22}^l&= s\frac\partial{\partial s}+(1/2-\mu_2).
\label{eq:El2}
\end{align}

{\bf \punct. Correspondence of differential operators.}  
We want to evaluate images of differential operators under the Fourier transform, i.e.
for a differential operator $D$ in $C_0^\infty\bigl(\GL(2,\R)\bigr)$ we wish to find
a transformation $\Theta(D)$ in $\cE$ such that
$$
K_{DF}=\Theta(D) K_F
$$
for all $F$.

\sm

Define the following shift operators 
 \begin{align*}
 V_1  K(t,s|\mu_1,\epsilon_1;\mu_2,\epsilon_2)= K(t,s|\mu_1+ 1,\epsilon_1+1;\mu_2,\epsilon_2);
 \\
 V_2  K(t,s|\mu_1,\epsilon_1;\mu_2,\epsilon_2)= K(t,s|\mu_1,\epsilon_1;\mu_2+ 1,\epsilon_2+ 1).
 \end{align*}

\sm

{\bf Theorem 1.} a) 
{\it
The Fourier-images of operators of multiplication
by functions $a$, $b$, $c$, $d$, $(ad-bc)^{-1}$ in $C_0^\infty\bigl(\GL(2,\R)\bigr)$
 are the following operators in $\cE$:} 
\begin{align}
 \label{eq:th1b}
a&\quad \longleftrightarrow\quad  V_1 - \frac t {\mu_1-\mu_2}\Bigl(\frac\partial {\partial t} V_1 +  \frac\partial {\partial s} V_2  \Bigr);
\\
b&\quad \longleftrightarrow\quad
s V_1 - t V_2 -\frac{st} {\mu_1-\mu_2} \Bigl(\frac\partial {\partial t} V_1 +  \frac\partial {\partial s} V_2  \Bigr);
\\
c& \quad \longleftrightarrow\quad 
\frac{1} {\mu_1-\mu_2} \Bigl(\frac\partial {\partial t} V_1 +  \frac\partial {\partial s} V_2  \Bigr);
\\
d &\quad \longleftrightarrow\quad V_2 + \frac{s} {\mu_1-\mu_2} \Bigl(\frac\partial {\partial t} V_1 +  \frac\partial {\partial s} V_2  \Bigr);
\\
(ad-bc)^{-1} &\quad \longleftrightarrow\quad V_1^{-1} V_2^{-1}.
 \label{eq:th1l}
\end{align}

b) {\it The operators $\frac\partial {\partial a}$,  $\frac\partial {\partial b}$,
$\frac\partial {\partial c}$, $\frac\partial {\partial d}$
correspond to:}
\begin{align}
 \label{eq:th2b}
 \frac\partial {\partial a} &\quad \longleftrightarrow\quad 
 (\tfrac32-\mu_1) V_1^{-1} + \frac s {\mu_1-\mu_2}\Bigl((\tfrac32-\mu_1)  \frac\partial {\partial s} V_1^{-1} + 
 (\tfrac32-\mu_2)  \frac\partial {\partial t} V_2^{-1}   \Bigr);
 \\
 \frac\partial {\partial b} &\quad \longleftrightarrow\quad - 
 \frac 1 {\mu_1-\mu_2}\Bigl((\tfrac32-\mu_1)  \frac\partial {\partial s} V_1^{-1} + 
 (\tfrac32-\mu_2)  \frac\partial {\partial t} V_2^{-1}   \Bigr);
 \\
 \frac\partial {\partial c} &\quad \longleftrightarrow\quad 
 (\tfrac32-\mu_1) tV_1^{-1} - (\tfrac32-\mu_2) sV_2^{-1}+ 
 \\ &
 \phantom{\frac\partial {\partial c} \quad \longleftrightarrow\quad\qquad }
 +
 \frac {st} {\mu_1-\mu_2}\Bigl((\tfrac32-\mu_1)  \frac\partial {\partial s} V_1^{-1} + 
 (\tfrac32-\mu_2)  \frac\partial {\partial t} V_2^{-1}   \Bigr);
\nonumber
 \\
  \frac\partial {\partial d} &\quad \longleftrightarrow\quad 
  (\tfrac 32-\mu_2) V_2^{-1} -  \frac {t} {\mu_1-\mu_2}\Bigl((\tfrac32-\mu_1)  \frac\partial {\partial s} V_1^{-1} + 
 (\tfrac32-\mu_2)  \frac\partial {\partial t} V_2^{-1}   \Bigr).
 \label{eq:th2l}
\end{align}

{\sc Proof.} 
A derivation of the formula for $(ad-bc)^{-1}$ is straightforward.
Other statements follow from the table of formulas \cite{Ner-jfaa}, Subsection 2.5.
The table contains a correspondence for 8 pencils of operators of type $A+\sigma B$,
where $\sigma$ is a parameter. Substituting $\tau=-1+\sigma$ and considering 
correspondences of constant terms and terms of the first order in $\tau$, we found 
the desired correspondences.

Arguments of \cite{Ner-jfaa} can be easily transformed to direct proofs of our formulas,
but this does   not essentially simplifies calculations.
\hfill $\square$

\sm

As an immediate corollary of Theorem 1 we get the following statement.

\sm

{\bf Theorem 2.} {\it Let $D$ be a  differential operator  on the group $\GL(2,\R)$ admitting a 
polynomial expression in $a$, $b$, $c$, $d$, $(ad-bc)^{-1}$,
$\frac\partial {\partial a}$,  $\frac\partial {\partial b}$,
$\frac\partial {\partial c}$, $\frac\partial {\partial d}$.
Then its Fourier-image 
is a finite sum of the form
$$
\Theta(D)=
\sum_{k,l} Q_{kl}\Bigl(t,s,\frac\partial {\partial t}, \frac\partial {\partial s}\Bigr) \, V_1^k V_2^l, 
$$
where $k$, $l\in\Z$, and $Q_{kl}(\cdot)$ are polynomial expressions
with rational coefficients depending on $\mu_1$, $\mu_2$
with poles at lines $\mu_1-\mu_2=m\in\Z$.}

\sm

{\sc Remark.}
Similar statements hold for $\GL(2,\C)$, in fact formulas for correspondence 
of differential and differential-difference operators in this case are
 same, see \cite{Ner-jfaa}. \hfill $\boxtimes$
 
 \sm
 
 {\bf \punct. The problem about overalgebra.} In \cite{Ner-imaginary} there was formulated the following question.
 
 \sm
 
 {\bf Problem.} {\it Let $G$ be  a semisimple Lie group, $\frg$ its Lie algebra,
 $H\subset G$ a subgroup,
 and $\rho$ a unitary representation of $G$. Assume that a restriction of $\rho$ to $H$ admits an
 explicit spectral decomposition. To transfer the action of the whole
 Lie algebra $\frg$ in the spectral decomposition.}
  
 \sm
 
 In \cite{Ner-imaginary} there was obtained an explicit solution of the problem
 for a tensor product of a highest weight representation of  $\SL(2,\R)$ and  the dual
  lowest weight representation%
  \footnote{Problem of decomposition of a  tensor product $\rho_1\otimes \rho_2$
  of unitary representations
  of a group $G$ is a special case of a restriction problem, namely we restrict a representation
  of $G\times G$ to the diagonal $G$.}. The overalgebra $\mathfrak{sl}(2,\R)\oplus \mathfrak{sl}(2,\R)$
 act by differential-difference operators, including second derivatives. Now there is a collection 
 of solved problems of this type, see \cite{Mol1}--\cite{Mol3}, \cite{Ner-Gl}, \cite{Ner-jfaa}.
 In all the cases we observe differential-difference operators, in all cases we have shifts in
 imaginary directions. Results of \cite{Ner-Gl}
show that orders of partial derivatives increases with growth of a rank of a group.
In my opinion, now there are reasons to hope that such problems are always solvable
(as far as we can explicitly solve the restriction problem).

Formulas for $\GL(2,\R)$ mentioned in the proof of Theorem 1 were obtained as 
a solution of the overalgebra problem
for restrictions of representations of degenerate principal series%
\footnote{For any classical Lie group $G=\GL(n,\R)$, $\GL(n,\C)$, $\GL(n,\H)$,
$\O(p,q)$, $\U(p,q)$, $\Sp(p,q)$,   $\Sp(2n,\R)$, $\Sp(2n,\C)$, $\O(n,\C)$, $\SO^*(2n)$
the regular representation
of $G\times G$ in $L^2(G)$ can be obtained as a restriction from a certain unitary representation
of a certain overgroup $\wt G\supset G\times G$, namely, $\wt G= \GL(2n,\R)$, $\GL(2n,\C)$, $\GL(2n,\H)$,
$\O(p+q,p+q)$, $\U(p+q,p+q)$, $\Sp(p+q,p+q)$,  $\Sp(4n,\R)$, $\Sp(4n,\C)$, $\O(2n,\C)$, $\SO^*(4n)$.
More generally, the representation
of a classical group $H$ in $L^2$ on a  (pseudo-Riemannian) symmetric space $H/L$ can be obtained as a restriction
of a certain unitary representation of a certain overgroup $\wt H$. For details and precise 
formulations, see \cite{Ner-uniform}.}
of $\GL(4,\R)$
to $\GL(2,\R)\times \GL(2,\R)$. 

\sm 

{\bf \punct. The Weil representation.} Consider the real symplectic group
$\Sp(2n,\R)$, i.e., the group of block $(n+n)\times(n+n)$ real matrices 
$g=\begin{pmatrix}
  a&b\\c&d
 \end{pmatrix}
 $ satisfying the condition
 $$
 g \begin{pmatrix} 0&1\\-1&0\end{pmatrix} g^t= \begin{pmatrix} 0&1\\-1&0\end{pmatrix},$$
 the symbol $^t$ denotes the transposition. Its Lie algebra
 $\mathfrak{sp}(2n)$
 consists of real matrices of the form
 \begin{equation}
 \begin{pmatrix}
                    \alpha&\beta\\ \gamma&-\alpha^t
                   \end{pmatrix}, \qquad \text{where $\beta=\beta^t$, $\gamma=\gamma^t$.}
  \label{eq:Lie}
 \end{equation}
The matrices of the form $\begin{pmatrix}
                           a&0\\0&a^{t-1}
                          \end{pmatrix}$
form a subgroup in $\Sp(2n,\R)$, it is isomorphic to $\GL(n,\R)$.   

We also consider the Heisenberg group $\Heis_{2n+1}$ consisting of block matrices 
of size $(1+n+1)$ of the form
$\begin{pmatrix}
  1&v&\lambda\\ 0&1&w\\0&0&1
 \end{pmatrix}.
$

These groups have standard representations in the space $L^2(\R^n)$,
we briefly describe them on the level of Lie algebras.
The  collection of differential operators
$$
x_1,\dots, x_n, \frac{\partial}{\partial x_1}, \dots, \frac{\partial}{\partial x_n},  1
$$
is closed with respect to the commutator and
form a representation of the Heisenberg Lie algebra $\mathfrak{heis_{2n+1}}$.

The collection of operators (they are quadratic expressions in operators of the Heisenberg algebra) 
$$
 x_k x_l,\qquad  x_k \frac{\partial}{\partial x_l} + \delta_{kl},
\qquad \frac{\partial^2}{\partial x_k \partial x_l}
$$
also is closed with respect to the commutator and
forms the symplectic Lie algebra $\mathfrak{sp}(2n)$. This representation  is called 
the {\it Weil representation} or {\it oscillator representation} 
(the symplectic group $\Sp(2n,\R)$ itself
acts by integral operators,  see, e.g., \cite{Ner-gauss}).

The operators  $x_k \frac{\partial}{\partial x_l} +\delta_{kl}$ generate a representation of the
Lie algebra $\mathfrak{gl}(n)\subset \mathfrak{sp}(2n)$; the corresponding representation
of the Lie group $\GL(n,\R)$ is the standard representation of $\GL(n,\R)$  in
$L^2(\R^n)$).

Joining $\mathfrak{heis}_{2n+1}$ and $\mathfrak{sp}(2n)$
we get a representation of a semidirect product $\mathfrak{sp}(2n)\ltimes \mathfrak{heis}_{2n+1}$ of Lie algebras,
and a representation of the corresponding group $\Sp(2n,\R)\ltimes \Heis_{2n+1}$.

Recall that ($\GL(n,\R)$, $\GL(n,\R)$) is an example of Howe dual pair (see \cite{Howe}). 
We do not need a definition of such pairs and only present minimal data necessary for our purposes.

Consider  a symplectic group $\Sp(2 n^2,\R)$ and its subgroup $\GL(n^2,\R)$.
We realize $\R^{n^2}$ as a space $\Mat(n,\R)$ of real matrices of $n\times n$. The group
$\GL(n,\R)\times \GL(n,\R)$ acts on $\Mat(n,\R)$
by 
$$
(g_1,g_2):\, X\mapsto g_1^{-1} X g_2, \qquad \text{where $g_1$, $g_2\in \GL(n,\R)$.}
$$
 Therefore, we have an homomorphism $i:\GL(n,\R)\times \GL(n,\R)\to \GL(n^2,\R)$.
 We  restrict the Weil representation of $\Sp(2n^2,\R)$ to the subgroup
 $\GL(n,\R)\times \GL(n,\R)$ (or, more precisely, we consider the composition of
 $i$ and the Weil representation).
  Clearly, this restriction is the natural representation
 of $\GL(n,\R)\times \GL(n,\R)$ in $L^2\bigl(\Mat(n,\R)\bigr)$. On the other hand
 the group $\GL(n,\R)$ is a dense open subset in $\Mat(n,\R)$
 and we can identify $L^2\bigl(\Mat(n,\R)\bigr)$ with the space
 $L^2$ on $\GL(n,\R)$ with respect to the Haar measure $\det(X)^{-n-1}\,dX$,
 the identification is given by the formula
 $$
 f\mapsto |\det(X)|^{(n+1)/2} f(X),
 $$
 the factor $\det(X)^{-(n+1)/2} $ is necessary to make the operator unitary.
 There arise a question: 
 
 \sm
 
{\it  To obtain the action of the Lie algebra 
$\mathfrak{sp}(2n^2)\ltimes \mathfrak{heis}_{2n^2+1}$ in the spectral decomposition of $L^2\bigl(\GL(n,\R)\bigr)$.}

\sm

Clearly, Theorem 1 gives a solution for $n=2$.

\sm

{\bf \punct. Action of $\mathfrak{sp}(8)$ in the space $\cE$.}
Return to the notation $a$, $b$, $c$, $d$ for coordinates on $\GL(2,\R)$ and $\Mat(2,\R)$.
The Heisenberg algebra acts on $\Mat(2,\R)$ by the operators
$$\text{$ia$, $ib$, $ic$, $id$, $\frac{\partial}{\partial a}$, $\frac{\partial}{\partial b}$, 
$\frac{\partial}{\partial c}$, $\frac{\partial}{\partial d}$.}$$
We identify $L^2(\Mat(2,\R))$ and $L^2$ on the group by
 $$f\begin{pmatrix}a&b\\c&d\end{pmatrix}
   \,\mapsto\,
 \begin{pmatrix}a&b\\c&d\end{pmatrix}
     |ad-bc|. $$ 
 The coordinate operators $ia$, $ib$, $ic$, $id$ under this identification
remain the same, the partial derivatives  transform
to
\begin{align*}
 \partial_a:= \frac{\partial}{\partial a}- \frac d{ad-bc};
 \qquad
  \partial_b:= \frac{\partial}{\partial b}+ \frac c{ad-bc};
  \\
   \partial_c:= \frac{\partial}{\partial c}+ \frac b{ad-bc};
   \qquad
    \partial_d:= \frac{\partial}{\partial d}- \frac a{ad-bc}.
\end{align*}
This gives us the action of the Heisenberg algebra $\mathfrak{heis}_9$.
The  operators in $\cE$ corresponding to $\partial_a$, $\partial_b$, etc.
can be easily evaluated. A straightforward calculation shows that they can be obtained from formulas
(\ref{eq:th2b})--(\ref{eq:th2l})
by a substitution $\mu_1\to \mu_1+1$, $\mu_2\to \mu_2+1$
(this is equivalent to a formal substitution $\frac 32\mapsto \frac 12$ to expressions (\ref{eq:th2b})--(\ref{eq:th2l})).

This gives us an action of the Heisenberg algebra $\mathfrak{heis}_9$ in $\cE$.
The symplectic Lie algebra
 $\mathfrak{sp}(8)$ acts by quadratic expressions in generators
of the Heisenberg algebra. The author does not know a nice general formula for 
all 36 generators of $\mathfrak{sp}(8)$, but for a given element  its expression can be
written. For instance,
\begin{multline*}
c^2 \quad \longleftrightarrow \quad 
\frac 1{(\mu_1-\mu_2)(\mu_1-\mu_2+1)} \frac {\partial^2}{\partial t^2} V_1^2
+ \\ +
 \frac 2{(\mu_1-\mu_2-1)(\mu_1-\mu_2+1)} \frac {\partial^2}{\partial t\,\partial s} V_1 V_2
 +
  \frac 1{(\mu_1-\mu_2)(\mu_1-\mu_2-1)} \frac {\partial^2}{\partial s^2}  V_2^2;
\end{multline*}


\begin{multline*}
 c\partial_b \quad \longleftrightarrow \quad \frac{\mu_2-\tfrac 12}{(\mu_1-\mu_2+1) (\mu_1-\mu_2) }
 \frac{\partial^2}{\partial t^2} V_1 V_2^{-1}
 +\\+
 \frac{\mu_1+\mu_2}{(\mu_1-\mu_2+1)(\mu_1-\mu_2-1)} \frac {\partial^2}{\partial t\,\partial s} 
 +\frac{\mu_1-\tfrac 12}{(\mu_1-\mu_2)(\mu_1-\mu_2-1)} \frac{\partial^2}{\partial s^2} V_1^{-1} V_2
\end{multline*}

\begin{multline*}
 \partial^2_b \quad \longleftrightarrow \quad 
 \frac{(\mu_2-\tfrac 32) (\mu_2-\tfrac 12)  }{(\mu_1-\mu_2+1) (\mu_1-\mu_2)} \frac{\partial^2}{\partial t^2}
 V_2^{-1}+
 \\+
 \frac{2 (\mu_2-\tfrac 12)(\mu_1-\tfrac 12)}{(\mu_1-\mu_2+1) (\mu_1-\mu_2-1)}
 \frac{\partial^2}{\partial t\, \partial s} V_1^{-1} V_2^{-1}
 +
 \frac{(\mu_1-\tfrac 32) (\mu_1-\tfrac 12) }{(\mu_1-\mu_2)(\mu_1-\mu_2-1)}
 \frac{\partial^2}{\partial t^2} V_1^{-2}.
\end{multline*}

 For some quadratic operators Fourier-images have a simpler form. Namely
for generators (\ref{eq:er1})--(\ref{eq:el2}) of the subalgebra $\mathfrak{gl}(2)\oplus \mathfrak{gl}(2)$
we have operators (\ref{eq:Er1})--(\ref{eq:El2}).
In particular,
\begin{equation}
a\frac \partial{\partial a}+ b\frac \partial{\partial b}+ c\frac \partial{\partial c}+
d\frac \partial{\partial d}  \quad \longleftrightarrow \quad  -\mu_1-\mu_2.
\label{eq:inv1}
\end{equation}
 The operator corresponding to
 $(ad-bc)$ is given by (\ref{eq:th1l}). Also
\begin{equation}
\partial_a\partial_d-\partial_b\partial_c \quad \longleftrightarrow \quad 
\bigl(\mu_1-\tfrac 12\bigr)\bigl(\mu_2-\tfrac12\bigr) V_1^{-1} V_2^{-1}.
\label{eq:inv2}
\end{equation}

{\sc Remark.} These simple expressions are related to a structure of the Lie algebra $\mathfrak{sp}(8)$
as a $\mathfrak{sl}(2)\oplus \mathfrak{sl}(2)$-module.
First, we notice that algebra of matrices (\ref{eq:Lie})
splits into a direct sum of  subalgebras consisting of matrices 
$\begin{pmatrix}
                                      \alpha&0\\0&-\alpha^t
                                     \end{pmatrix}$,
 $\begin{pmatrix}
                                      0&\beta\\0&0
                                     \end{pmatrix}$,
 $\begin{pmatrix}
                                      0&0\\\gamma&0
                                     \end{pmatrix}$
                                     respectively.
 These subalgebras are invariant with respect to the adjoint action of                                   
$\mathfrak{sl}(2)\oplus \mathfrak{sl}(2)$.
 Denote by $W_k$ a  $k$-dimensional irreducible
$\mathfrak{sl}(2)$-module. 
The action of  $\mathfrak{sl}(2)\oplus \mathfrak{sl}(2)$
on symmetric matrices $\beta$ splits as 
$$(V_3\otimes V_3)\oplus (V_1\otimes V_1),$$
the same structure has the space of symmetric matrices $\gamma$.
The space of matrices $\alpha$ is
$$
(V_2\oplus V_2)\otimes (V_2\oplus V_2)= (V_3\otimes V_3)\oplus (V_3\otimes V_1)
\oplus (V_1\otimes V_3)\oplus (V_1\otimes V_1).
$$
For 3 invariants (elements of the spaces $V_1\otimes V_1$) we have formulas 
(\ref{eq:th1l}),
(\ref{eq:inv1}), (\ref{eq:inv2}). Two summands $V_3\otimes V_1$, $V_1\otimes V_3$
correspond to the subalgebra $\mathfrak{sl}(2)\oplus \mathfrak{sl}(2)$.
Three long formulas above correspond to representatives of three modules $V_3\otimes V_3$.
\hfill $\boxtimes$

\sm

{\bf  \punct. A toy example: the action of $\SL(2,\R)$ in functions on plane.}
Consider the natural action of the group $\SL(2,\R)$
$$
T\begin{pmatrix}
  a&b\\c&d
 \end{pmatrix}\phi(x,y)=\phi(ax+cy, bx+dy)
$$
on the space $C_0^\infty(\R^2)$  of smooth compactly supported functions on the plane 
$\R^2$.
For  $\phi(x,y)\in C_0^\infty(\R^2)$  we assign the function
\begin{equation}
F(u, \sigma,\epsilon)=J\phi(u, \sigma,\epsilon)  :=\int_\R \phi(t,tu) \, t^{-\sigma\kos\epsilon}\, dt,
\label{eq:Fphi}
\end{equation}
 where $u\in \R$, $\sigma\in\C$, $\epsilon=0$, $1$. 
 For 
 $\Re\sigma<1$ the integral converges and is holomorphic in $\sigma$,
 this function admits a meromorphic continuation to arbitrary $\sigma$   with possible poles at $\sigma=1$, 2, 3, \dots.
 Denote by $\cM$ the space of all functions $F$ that can be obtained in this way. 
The action of $\SL(2,\R)$ transfers to the action
$$
R\begin{pmatrix}
  a&b\\c&d
 \end{pmatrix} F(u,\sigma,\epsilon)=F\Bigl( \frac{b+u d}{a+uc}, \,\sigma,\epsilon \Bigr) (a+uc)^{-1+\sigma\kos \epsilon}.
$$

Next, define the inner product on the space $\cM$  by
$$
\la F_1,F_2\ra=\sum_{\epsilon=1,2}\int_\R\int_\R F_1(u,i\tau,\epsilon) \ov{F_2(u,i\tau,\epsilon)}
du\,d\tau.
$$
Passing to the completion,
we get a Hilbert space $\cH$ (a sum of two copies of $L^2(\R^\times\R)$).
The operator $J$ extends to a unitary operator $J:L^2(\R^2)\to \cH$
(see, e.g., \cite{GGP}, \S I.3), and this gives us a spectral decomposition
of the representation $T$.

As we mentioned above, the representation $T$ is a restriction of the Weil
representation of $\Sp(4,\R)\ltimes \mathrm{Heis}_5$ to the subgroup
$\SL(2,\R)$. Our previous considerations suggest that we can 
transfer the action of the Lie algebra
 $\mathfrak{sp}(4)\ltimes \mathfrak{heis}_5$ to the space $\cM$.
 
 Define a shift operator
 $$
 V F(u,\sigma,\epsilon)=F(u,\sigma+1,\epsilon+1).
 $$

We have the following correspondence of  operators:
\begin{align}
 x &\quad \longleftrightarrow \quad V^{-1};
  \label{eq:x}
 \\
 y &\quad \longleftrightarrow \quad u V^{-1};
 \label{eq:y}
 \\
 \frac\partial{\partial x}& \quad \longleftrightarrow \quad  \sigma V- u \frac\partial{\partial u} V;
  \label{eq:dx}
 \\
 \frac\partial{\partial y}& \quad \longleftrightarrow \quad u \frac\partial{\partial u} V.
  \label{eq:dy}
\end{align}
 
 Let us verify (\ref{eq:dy}). By (\ref{eq:Fphi}), we get
 $$
\frac \partial{\partial u} F(u, \sigma+1,\epsilon+1):=
\int_\R t \frac{\partial \phi}{\partial u} (t,tu) \, t^{-\sigma+1\kos\epsilon+1}\, dt
= \int_\R  \frac{\partial \phi}{\partial u} \phi(t,tu) \, t^{-\sigma\kos\epsilon}\, dt,
 $$
this is the required result. A derivation of formulas (\ref{eq:x}) and (\ref{eq:y})
are straightforward.
To verify (\ref{eq:dx}) we write
\begin{equation}
\int_\R \frac d{dt} \bigl(\phi(t,tu)\bigr) \, t^{-\sigma\kos\epsilon}\, dt=
\int_\R  \frac{\partial \phi}{\partial x} (t,t u)  \, t^{-\sigma\kos\epsilon}\, dt
+ \int_\R  (tu) \frac{\partial \phi}{\partial y} (t,t u)  \, t\cdot t^{-\sigma\kos\epsilon}\, dt.
\label{eq:last}
\end{equation}
For $\Re\sigma<0$,
integrating the left hand side by parts we get $\sigma F(u,\sigma+1,\epsilon+1)$,
by analytic continuation we extend this to arbitrary $\sigma\in \C$.
The second term in the right hand side is 
$$
\int_\R   \Bigl(y\frac{\partial  \phi(x, y)}{\partial y}\Bigr)\Bigr|_{x=t,y=tu} 
\, t\cdot t^{-\sigma+1\kos\epsilon+1}\, dt. 
$$
By (\ref{eq:y}) and (\ref{eq:dy}) it equals $u\frac\partial{\partial u} F(u,\sigma+1,\epsilon+1)$.
The first term in the right hand side of (\ref{eq:last}) is the $J$-image of $\frac\partial{\partial x}$.

So we have a possibility to evaluate  the image of any differential operator, which is
polynomial in
$x$, $y$, $\frac{\partial}{\partial x}$, $\frac{\partial }{\partial y}$.
In particular, we get formulas for the action of the symplectic Lie algebra $\mathfrak{sp}(4)$.

 \tt
 
 \noindent
 Math. Dept., University of Vienna; \\
 Institute for Theoretical and Experimental Physics (Moscow); \\
 MechMath Dept., Moscow State University;\\
 Institute for Information Transmission Problems;\\
 yurii.neretin@univie.ac.at;
 URL: http://mat.univie.ac.at/$\sim$neretin/
 

\end{document}